\documentclass{amsart}
\usepackage{verbatim}
\usepackage{paralist}
\usepackage{leftidx}
\usepackage{dsfont}
\usepackage{amsthm}
\usepackage{amsmath}
\usepackage{amssymb}
\usepackage{amscd}
\usepackage{graphicx}
\usepackage{setspace}
\usepackage[all]{xy}
\usepackage{mathrsfs} 
\usepackage{hyperref}
\usepackage{mathtools}
\usepackage[font=small,labelfont=bf]{caption} 
\title{On the Faltings height of the curve $y^2=x^n-1$}
\author{Robert Wilms}
\address{Universit\'e de Caen Normandie, CNRS, LMNO UMR 6139, F-14000 Caen, France}
\email{\href{mailto:robert.wilms87@gmail.com}{robert.wilms87@gmail.com}}
\date{\today}
\thanks{The author gratefully acknowledges funding from the \emph{Région Normandie} under the program \emph{``Normandie Recherche -- Objectif Labels d’excellence 2025''} (Project \textbf{MERCI})}
\numberwithin{equation}{section}
\newtheorem{Def}{Definition}
\numberwithin{Def}{section}

\newtheorem{Rem}[Def]{Remark}
\newtheorem{Lem}[Def]{Lemma}

\newtheorem{Cor}[Def]{Corollary}

\newtheorem{Thm}[Def]{Theorem}

\begin{document}
	\begin{abstract}
		We compute the stable Faltings height of the hyperelliptic curve
		$X_n\colon y^2=x^{n}-1$ for every odd integer $n\ge 3$ in terms of special values of Euler's gamma function.
		In particular, we prove the bounds
		$$-0.975n< h_{\mathrm{Fal}}(X_n)-\tfrac{n}{8}\log n<\tfrac{9}{64}n\log\log n-0.263n.$$
		As an application, we bound the Faltings height of any abelian variety with complex multiplication by the canonical CM-type of the $n$-th cyclotomic field by $\frac{n}{8}\log n+\frac{9}{64}n\log\log n-0.136n$.
	\end{abstract}
	
	\maketitle
	\section{Introduction}
	In his proof of the Mordell conjecture \cite{Fal83}, Faltings introduced a height function on abelian varieties
	with the Northcott property and good behavior under isogenies. This played a key role in his proof of the Shafarevich conjecture, namely the finiteness of isomorphism classes
	of abelian varieties over a fixed number field with fixed dimension, fixed locus of bad reduction, and fixed polarization degree.
	Since then, many aspects of the Faltings height have been studied.
	For instance, Colmez \cite{Col93} investigated the case of complex multiplication and conjectured a formula in terms
	of special values of $L$-functions.
	Bost \cite{Bos96} proved an explicit lower bound, and Burgos Gil, Menares, and Rivera-Letelier \cite{BMR18}
	studied the essential minimum for elliptic curves.
	The generalized Szpiro conjecture predicts an upper bound for the Faltings height in terms of the norm of the conductor
	and the discriminant of the field of definition. We refer to \cite{GKM19} for a discussion on this conjecture.
	
	The aim of this note is to compute the Faltings height of the family of hyperelliptic curves defined by the equation
	\begin{align}\label{equ_curve-Xn}
		X_n\colon\qquad y^2=x^{n}-1		
	\end{align}
	for all odd $n\ge 3$. In general, one defines the Faltings height of a curve $X$ by the Faltings height of its Jacobian, that is, $h_{\mathrm{Fal}}(X)=h_{\mathrm{Fal}}(\mathrm{Jac}(X))$. We write $\mathrm{ord}_p(n)$ for the exponent of $p$ in the prime factorization of $n$ and $\Gamma$ denotes Euler's gamma function. The following theorem gives an explicit expression for the Faltings height of $X_n$ in terms of values of the gamma function.
	\begin{Thm}\label{thm_faltings-height-hyperelliptic}
		Let $X_n$ be the hyperelliptic curve defined by the equation $y^2=x^n-1$ for any odd integer $n\ge 3$. We write $g=\frac{n-1}{2}$ for its genus. The stable Faltings height of $X_n$ is given by
		\begin{align*}
			h_{\mathrm{Fal}}(X_n)=&\tfrac{n}{8}\sum_{\substack{p\text{ prime}\\ p\mid n}}\frac{p^{2\mathrm{ord}_{p}(n)}-1}{p^{2\mathrm{ord}_{p}(n)-1}(p^2-1)}\log p+\tfrac{n}{8}\log n-\tfrac{g}{2}\log \pi+\tfrac{g}{2n}\log 2\\
			&-\sum_{j=1}^g\log\frac{\Gamma\left(\frac{2j-1}{2n}\right)}{\Gamma\left(\frac{1}{2}+\frac{2j-1}{2n}\right)}.
		\end{align*}
	\end{Thm}
	Note that the \emph{stable} Faltings height is the Faltings height computed over a base field, where the curve has semistable reduction. We will prove the theorem by an explicit calculation of the vanishing orders at the finite places and the norm at the infinite places of a wedge product of a canonical basis of differential forms in $H^0(X_n,\Omega_{X_n}^1)$. For the computation of the vanishing orders at the finite places, we will use the theory of cluster pictures introduced by Dokchitser, Dokchitser, Maistret, and Morgan in \cite{DDMM23} as well as an explicit result by Kunzweiler \cite{Kun20}. The norm at the infinite places is in general given by a determinant of a $g\times g$ matrix, which we can explicitly compute by reduction to a Vandermonde determinant in this special situation. Note that for $n=5$ the result has already been shown by Bost--Mestre--Moret--Bailly \cite{BMM90} and our computation is inspired by their work.
	
	\begin{figure}[t!]
		\centering
		\includegraphics[width=0.96\textwidth]{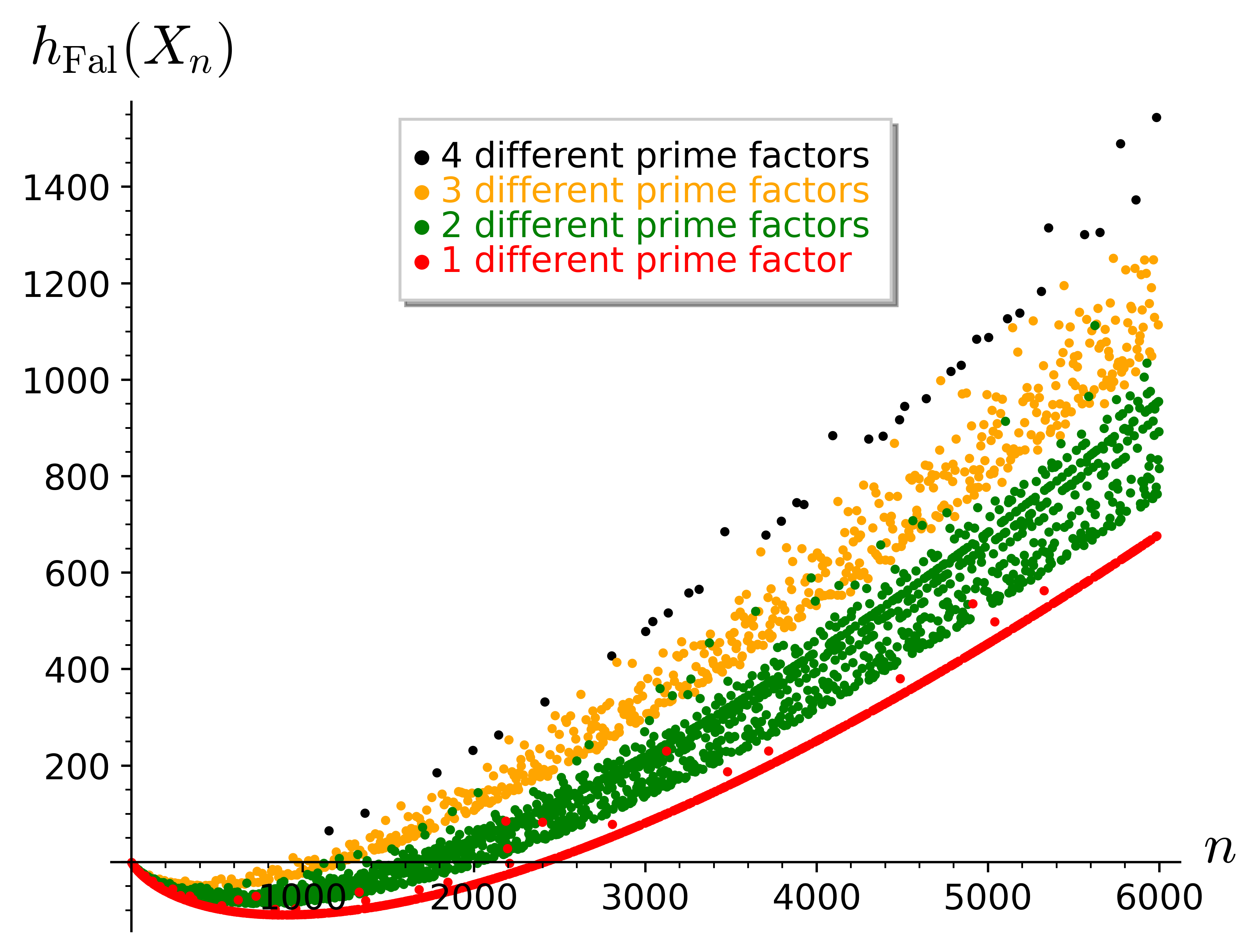}
		\caption{Plot of the points $(n,h_{\mathrm{Fal}}(X_n))$ for all odd $3\le n\le 6001$.}
		\label{fig_faltings-height-hyperelliptic}
	\end{figure}
	Theorem \ref{thm_faltings-height-hyperelliptic} allows us to plot the values of the Faltings height for $X_n$. Figure \ref{fig_faltings-height-hyperelliptic} shows the values of the Faltings height of $X_n$ for all odd integers $n$ from $3$ to $6001$. As a corollary of Theorem \ref{thm_faltings-height-hyperelliptic}, we give explicit lower and upper bounds for the stable Faltings height $h_{\mathrm{Fal}}(X_n)$ and determine its asymptotic behavior with respect to $n$.	
	\begin{Cor}\label{cor_bounds}
		For any odd integer $n\ge 3$, the stable Faltings height $h_{\mathrm{Fal}}(X_n)$ of the curve $X_n$ satisfies the bounds
		$$-0.975n< h_{\mathrm{Fal}}(X_n)-\tfrac{n}{8}\log n<\tfrac{9}{64}n\log\log n-0.263n.$$
		In particular, we have the asymptotic behavior
		$$h_{\mathrm{Fal}}(X_n)=\tfrac{1}{8}n\log n +O(n\log\log n).$$
	\end{Cor}
	In particular, the asymptotically leading term of $h_{\mathrm{Fal}}(X_n)$ does not depend on the prime factorization of $n$. In the proof we approximate the sum of the gamma functions in the expression of Theorem \ref{thm_faltings-height-hyperelliptic} by an integral and we bound the sum over the prime factors of $n$ using results on the distribution of primes by Rosser and Schoenfeld \cite{RS62}.
	
	Next, we discuss an application to abelian varieties with complex multiplication by a cyclotomic field. Let $n\ge 3$ be an odd integer and $\zeta_n\in\mathbb{C}$ a primitive $n$-th root of unity. The number field $\mathbb{Q}(\zeta_n)$ is a \emph{CM field}, that is a totally imaginary quadratic extension of a totally real number field. A \emph{CM type} $\Sigma$ of $\mathbb{Q}(\zeta_n)$ is a subset $\Sigma\subseteq \mathrm{Hom}(\mathbb{Q}(\zeta_n),\mathbb{C})$ such that $\Sigma\sqcup (\Sigma\circ c)=\mathrm{Hom}(\mathbb{Q}(\zeta_n),\mathbb{C})$, where $c\colon \mathbb{Q}(\zeta_n)\to\mathbb{Q}(\zeta_n)$ denotes the complex conjugation.
	Up to translation by a field automorphism in $\mathrm{Gal}(\mathbb{Q}(\zeta_n)/\mathbb{Q})$, we get a \emph{canonical CM type} 
	$$\Sigma_n=\{\sigma\in \mathrm{Hom}(\mathbb{Q}(\zeta_n),\mathbb{C})~|~\mathrm{Im}\,\sigma(\zeta_n)>0\}$$
	of $\mathbb{Q}(\zeta_n)$.
	
	We say that an abelian variety $A$ defined over $\mathbb{C}$ has \emph{complex multiplication by} $(\mathbb{Q}(\zeta_n),\Sigma)$ if there is a $\mathbb{Q}$-algebra embedding
	$$\iota\colon \mathbb{Q}(\zeta_n)\xhookrightarrow{} \mathrm{End}^0(A)=\mathrm{End}(A)\otimes_{\mathbb{Z}}\mathbb{Q},$$
	$A$ has dimension $\dim A=\frac{[\mathbb{Q}(\zeta_n):\mathbb{Q}]}{2}=\frac{\varphi(n)}{2}$, and the action of $\mathbb{Q}(\zeta_n)$ on the tangent space $T_0 A$ induced by $\iota$ yields an isomorphism
	$$T_0 A\cong\bigoplus_{\tau\in\Sigma}\mathbb{C}_{\tau}$$
	of $\mathbb{Q}(\zeta_n)\otimes_{\mathbb{Q}}\mathbb{C}$-modules, where $\mathbb{C}_{\tau}$ is a one-dimensional $\mathbb{C}$-linear space, on which $a\in \mathbb{Q}(\zeta_n)$ acts via $\tau(a)$. By a theorem of Colmez \cite[Theorem 0.3 (ii)]{Col93}, the stable Faltings height $h_{\mathrm{Fal}}(A)$ of an abelian variety $A$ with complex multiplication only depends on the corresponding CM type. We will show that for every prime $p\ge3$ the Jacobian $\mathrm{Jac}(X_p)$ has complex multiplication by $(\mathbb{Q}(\zeta_p),\Sigma_p)$, see Remark \ref{rem_CM}. Thus, Theorem \ref{thm_faltings-height-hyperelliptic} computes the stable Faltings height for these CM types. For general $n$, we will deduce the following upper bound for the Faltings height corresponding to the CM type $(\mathbb{Q}(\zeta_n),\Sigma_n)$ as an application of Theorem \ref{thm_faltings-height-hyperelliptic}.
	\begin{Cor}\label{cor_faltings-height-cyclotomic}
		Let $n\ge 3$ be an odd integer and let $A_n$ be an abelian variety having complex multiplication by the canonical CM type $(\mathbb{Q}(\zeta_n),\Sigma_n)$. The stable Faltings height of $A_n$ is bounded from above by
		\begin{align*}
		h_{\mathrm{Fal}}(A_n)< \tfrac{n}{8}\log n+\tfrac{9}{64}n\log\log n-0.136n. 
		\end{align*}
	\end{Cor}
	The corollary follows from Theorem \ref{thm_faltings-height-hyperelliptic} since we can establish $A_n$ as an abelian subvariety of $\mathrm{Jac}(X_n)$. It was proven by Rémond \cite{Rem22} that the Faltings height of an abelian subvariety is bounded by the Faltings height of the abelian variety up to an explicit error term.
	
	\subsubsection*{Outline}
	In Section \ref{sec_Faltings-height} we recall the definition of the stable Faltings height and we introduce a canonical section
	used to compute it.
	Section \ref{sec_cluster} reviews the cluster picture formalism of Dokchitser--Dokchitser--Maistret--Morgan and determines the
	clusters for $X_n$.
	In Section \ref{sec_finite} we compute the contributions of the finite places to the Faltings height.
	Sections \ref{sec_period-matrix} and \ref{sec_infinite} treat the archimedean contribution via an explicit computation of the period matrix.
	In Section \ref{sec_proof-main-thm} we assemble these ingredients to prove Theorem \ref{thm_faltings-height-hyperelliptic}.
	Section \ref{sec_bound} provides bounds for the relevant gamma products, which are then used in Section \ref{sec_proof-corollary} to
	prove Corollary \ref{cor_bounds}.
	Finally, Section \ref{sec_cm} proves Corollary \ref{cor_faltings-height-cyclotomic}.
	
	\section{The stable Faltings height}\label{sec_Faltings-height}
	We give the definition of the stable Faltings height of a hyperelliptic curve.	
	Let $K$ be a number field and $X$ a hyperelliptic curve of genus $g$ defined over $K$. This means that $X$ is a smooth projective curve given by the closure of the zero set of an equation
	\begin{align}\label{equ_polynomial-hyperelliptic}
	y^2=f(x)
	\end{align}
	in a weighted projective plane. Here, $f(x)$ is a monic separable polynomial of degree $2g+1$ or $2g+2$. After a coordinate change, we can and will assume that $\deg f=2g+2$. Let $\mathcal{X}\to B$ be the minimal regular model of $X$ over the spectrum $B=\mathrm{Spec}(\mathcal{O}_K)$ of the integers $\mathcal{O}_K$ of $K$. Replacing $K$ by some finite field extension, we may assume that $\mathcal{X}$ has everywhere semistable reduction. For every embedding $\sigma\colon K\to\mathbb{C}$ we write $X_\sigma$ for the hyperelliptic Riemann surface obtained by the base change induced by $\sigma$. We write $\omega_{\mathcal{X}/B}$ for the relative dualizing line bundle of $\mathcal{X}$ over $B$. The stable Faltings height is given by 
	$$h_{\mathrm{Fal}}(X)=\frac{1}{[K\colon \mathbb{Q}]}\widehat{\deg}\det H^0(\mathcal{X},\omega_{\mathcal{X}/B}),$$
	where $\widehat{\deg}$ denotes the Arakelov degree of the metrized line bundle $\det H^0(\mathcal{X},\omega_{\mathcal{X}/B})$ on $B$ and for every embedding $\sigma\colon K\to\mathbb{C}$ the line $\det H^0(X_{\sigma},\Omega_{X_{\sigma}}^1)$ is equipped with its Faltings metric $\|\cdot\|_\sigma$ induced by the inner product
	$$\langle \eta,\eta'\rangle =\frac{i}{2}\int_{X_{\sigma}}\eta\wedge \overline{\eta'}$$
	on $H^0(X_{\sigma},\Omega_{X_{\sigma}}^1)$. If $s$ is a non-zero rational section of $\det H^0(\mathcal{X},\omega_{\mathcal{X}/B})$, we have
	\begin{align}\label{equ_heightbysection}
	h_{\mathrm{Fal}}(X)=\frac{1}{[K\colon \mathbb{Q}]}\sum_{\mathfrak{p}\in |B|}\mathrm{ord}_\mathfrak{p}(s)\log N(\mathfrak{p})-\frac{1}{[K\colon \mathbb{Q}]}\sum_{\sigma\colon K\to\mathbb{C}}\log\|s\|_{\sigma},
	\end{align}
	where $|B|$ denotes the closed points of $B$ and for every $\mathfrak{p}\in |B|$ we write $N(\mathfrak{p})$ for the norm of the corresponding prime in $\mathcal{O}_K$.
	
	We can define a canonical section of the line bundle $\det H^0(X,\Omega_X^1)^{\otimes (8g+4)}$ by
	$$\Lambda=\frac{\mathrm{Disc}(f)^g}{2^{4g(g+1)}}\left(\frac{dx}{y}\wedge\frac{xdx}{y}\wedge\dots\wedge\frac{x^{g-1}dx}{y}\right)^{\otimes (8g+4)},$$
	where $\mathrm{Disc}(f)$ denotes the discriminant of the polynomial $f$. This section is independent of the choice of $f$ as shown in \cite[Proposition 2.2]{Kau99}. It naturally extends to a rational section on $\det H^0(\mathcal{X}, \omega_{\mathcal{X}/B})^{\otimes (8g+4)}$. Thus by the formula for the Arakelov degree as in Equation (\ref{equ_heightbysection}), we only have to compute $\mathrm{ord}_{\mathfrak{p}}(\Lambda)$ for every $\mathfrak{p}\in |B|$ and $\|\Lambda\|_\sigma$ for every embedding $\sigma\colon K\to\mathbb{C}$ to compute the Faltings height of $X$.
	
	\section{Cluster pictures}\label{sec_cluster}
	We recall the theory of cluster pictures introduced by Dokchitser, Dokchitser, Maistret, and Morgan in \cite{DDMM23}. This will allow us to compute the contributions of the finite places to the Faltings height in the next section. We assume that $\zeta_n\in K$.
	For any prime $\mathfrak{p}\subseteq \mathcal{O}_K$ we write $K_\mathfrak{p}$ for the completion of $K$ with respect to $\mathfrak{p}$ and $\mathcal{O}_{K_{\mathfrak{p}}}$ for the ring of integers of $K_{\mathfrak{p}}$. We write $\mathcal{X}_\mathfrak{p}$ for the base change of $\mathcal{X}$ induced by the embedding $\mathcal{O}_K\to\mathcal{O}_{K_{\mathfrak{p}}}$ and we set $B_\mathfrak{p}=\mathrm{Spec}(\mathcal{O}_{K_{\mathfrak{p}}})$.
	
	Now let $\mathfrak{p}\subseteq \mathcal{O}_{K}$ be any non-zero prime ideal with $2\notin\mathfrak{p}$. We denote by $v$ the valuation on $K_{\mathfrak{p}}$ normalized by $v(p)=1$ for the unique prime number $p\in\mathfrak{p}\cap\mathbb{Z}$. In particular, if $\pi_{\mathfrak{p}}\in\mathfrak{p}\mathcal{O}_{K_{\mathfrak{p}}}$ is a uniformizer, we have $\mathrm{ord}_{\mathfrak{p}}(\pi_{\mathfrak{p}})=1$ while $v(\pi_\mathfrak{p})=1/e_{\mathfrak{p}}$, where $e_{\mathfrak{p}}$ denotes the ramification index of $\mathfrak{p}$ over $p$. In general, we have $\mathrm{ord}_{\mathfrak{p}}(x)=e_{\mathfrak{p}}v(x)$ for all $x\in K_{\mathfrak{p}}$.
	We write $\mathfrak{R}$ for the set of roots of $f$ in an algebraic closure $\overline{K_{\mathfrak{p}}}$.
	By a \emph{cluster} we mean any non-empty subset $\mathfrak{s}\subseteq \mathfrak{R}$ such that there exists a disc $D=\{x\in\overline{K_{\mathfrak{p}}}~|~v(x-z)\ge d\}$ for some $d\in\mathbb{Q}$ and $z\in\overline{K_{\mathfrak{p}}}$ such that $\mathfrak{s}=D\cap \mathfrak{R}$. We call $\mathfrak{s}$ \emph{proper} if $\#\mathfrak{s}>1$. The \emph{depth} of a proper cluster $\mathfrak{s}$ is defined by $d_{\mathfrak{s}}=\min_{r,r'\in \mathfrak{s}}v(r-r')$. The \emph{relative depth} of a proper cluster $\mathfrak{s}\neq \mathfrak{R}$ is defined by $\delta_{\mathfrak{s}}=d_{\mathfrak{s}}-d_{\mathfrak{s}'}$ where $\mathfrak{s}'$ is the smallest cluster satisfying $\mathfrak{s}\subsetneq \mathfrak{s}'$.
	As an example, let us study the cluster picture of the curve $X_n$ given by Equation (\ref{equ_curve-Xn}) for any odd integer $n\ge 3$. 
	Let us first remark that the curve $X_n$ can also be defined by two natural alternative equations.
		\begin{Rem}\label{rem_coordinate-change}
		By a coordinate change $t=\frac{1}{x}$ and $z=\frac{y}{x^{(n+1)/2}}$ we see that $X_n$ is isomorphic over $\mathbb{Q}$ to the curve defined by the equation
		$$z^2=t(1-t^{n}).$$
		Likewise, by a coordinate change $u=(-4)^{-1/n}x$ and $v=\frac{y-i}{2i}$ we see that $X_n$ is isomorphic over $\mathbb{Q}(i,(-4)^{1/n})$ to the curve defined by the equation
		$$v^2+v=u^{n}.$$
	\end{Rem}
	Thus, we can choose the polynomial in (\ref{equ_polynomial-hyperelliptic}) to be $f_n(x)=x(1-x^{n})$. Note that $\mathfrak{R}=\{0,1,\zeta_n,\zeta_n^2,\dots,\zeta_n^{n-1}\}$, where $\zeta_n$ denotes a primitive $n$-th root of unity. The following lemma computes the clusters of $f_n$.
	\begin{Lem}\label{lem_cluster}
		Let $p$ be the unique prime number with $p\in\mathfrak{p}\cap\mathbb{Z}$.
		The clusters corresponding to $f_n$ are exactly given by $\{0\}$, $\mathfrak{R}$ and
		$$\mathfrak{s}_{a,b}=\left\{\left.\zeta_n^{\frac{kn}{p^a}+b}~\right|~k\in\mathbb{Z}/p^a\mathbb{Z}\right\}$$
		for $0\le a\le \mathrm{ord}_p(n)$ and $b\in \mathbb{Z}$. Moreover, we have $d_{\mathfrak{R}}=0$ and
		$$\delta_{\mathfrak{s}_{a,b}}=\begin{cases} \frac{1}{p^a} & \text{if } 1\le a<\mathrm{ord}_{p}(n),\\ \frac{1}{p^{a-1}(p-1)}&\text{if } 1\le a=\mathrm{ord}_{p}(n).\end{cases}$$
	\end{Lem}
	Note that the clusters $\mathfrak{s}_{a,b}$ only depend on the class of $b$ in $\mathbb{Z}/\frac{n}{p^a}\mathbb{Z}$.
	For the proof of the lemma we will need the following standard result on the valuations of differences of roots of unity.
	\begin{Lem}\label{lem_cyclotomic-valuation}
		Let $l\in\mathbb{Z}$ be arbitrary. It holds
		\begin{align*}
			v(1-\zeta_n^l)=\begin{cases} \infty &\text{if } n\mid l,\\ \frac{1}{p^{a-1}(p-1)}& \text{if } \frac{n}{\gcd(n,l)}=p^a \text{ for some } a\ge 1,\\ 0&\text{else}.\end{cases}
		\end{align*}
	\end{Lem}
	\begin{proof}
		We write $k=\frac{n}{\gcd(n,l)}$. In particular, $\zeta_n^l$ is a primitive $k$-th root of unity, which we denote by $\xi_k$. We have $k=1$ if and only if $n\mid l$ and if and only if $$v(1-\zeta_n^l)=v(0)=\infty.$$ Thus, we can assume $k\ge 2$. Let $\Phi_k(X)=\prod_{\substack{0< j< k\\ \gcd(j,k)=1}} (X-\xi_k^j)$ be the $k$-th cyclotomic polynomial. A standard result gives $$\prod_{\substack{0< j< k\\ \gcd(j,k)=1}} (1-\xi_k^j)=\Phi_k(1)=\begin{cases} q &\text{if } k=q^a \text{ for a prime } q \text{ and } a\in\mathbb{Z}_{\ge 1},\\ 1 & \text{else.}\end{cases}$$
		If $k$ is not a power of $p$, we conclude that $1-\zeta_n^l=1-\xi_k$ is a unit in $\mathcal{O}_{K_{\mathfrak{p}}}$ and hence $v(1-\zeta_n^l)=0$.
		If $k=p^a$ for some $a\in\mathbb{Z}_{\ge 1}$, then it follows from \cite[Lemma 10.1]{Neu99} that $p=u\cdot (1-\xi_k)^{\varphi(k)}$ for some unit $u\in \mathcal{O}_{K}^\times$. Since $\varphi(k)=p^{a-1}(p-1)$, we conclude $v(1-\zeta_n^l)=\frac{1}{p^{a-1}(p-1)}$. This completes the proof of the lemma.
	\end{proof}
	Now we give the proof of Lemma \ref{lem_cluster}.
	\begin{proof}[Proof of Lemma \ref{lem_cluster}]
		We first show that an arbitrary cluster $\mathfrak{s}$ of $f_n$ is equal to $\{0\}$, $\mathfrak{R}$, or $\mathfrak{s}_{a,b}$ for $0\le a\le \mathrm{ord}_p(n)$ and $b\in\mathbb{Z}$.
		Every singleton subset $\{r\}\subseteq \mathfrak{R}$ clearly defines a cluster and these are given by $\{0\}$ and $\mathfrak{s}_{0,b}$ for $b\in\{0,1,\dots,n-1\}$. Thus, we may assume that $\#\mathfrak{s}\ge 2$. Let $z\in \overline{K_{\mathfrak{p}}}$ and $d\in\mathbb{Q}$ be such that $$\mathfrak{s}=\{x\in \mathfrak{R}~|~v(x-z)\ge d\}.$$
		By the properties of ultrametrics, we can assume $z$ to be any point in $\mathfrak{s}$. By the assumption $\#\mathfrak{s}\ge 2$, we can moreover assume that $z=r_0$ for a root of unity $r_0\in\mathfrak{s}$. In particular, $v(z)=0$.		
		Let us first consider the case $0\in\mathfrak{s}$. This implies $v(z)\ge d$. Since all elements in $\mathfrak{s}$ have non-negative valuations, we get
		$$v(r-z)\ge \min\{v(r),v(z)\}=0\ge d$$
		for all $r\in\mathfrak{R}$ and hence $\mathfrak{s}=\mathfrak{R}$.
		
		Next, we consider the case $0\notin\mathfrak{s}$. By our assumptions, we have $z=\zeta_n^b$ for some $b\in\mathbb{Z}$. Every $r\in\mathfrak{R}\setminus\{0\}$ can be expressed in the form $r=\zeta_n^{b+t}$ for some $t\in\mathbb{Z}$ and hence, we get
		\begin{align}\label{equ_valuation-factor}
		v(r-z)=v(\zeta_n^{b+t}-\zeta_n^b)=v(\zeta_n^t-1).
		\end{align}
		In particular, the cluster $\mathfrak{s}$ can be expressed by
		$$\mathfrak{s}=\{\zeta_n^{b+t}\in\mathfrak{R}~|~v(\zeta_n^t-1)\ge d\}.$$
		We set
		$$a=\begin{cases} 0& \text{if } \frac{1}{p-1}<d,\\ a' & \text{if } \frac{1}{p^{a'}(p-1)}<d\le \frac{1}{p^{a'-1}(p-1)} \text{ for } a'\in\{1,\dots,\mathrm{ord}_p(n)\},\\ \mathrm{ord}_p(n)&\text{else.}\end{cases}$$
		It follows from Lemma \ref{lem_cyclotomic-valuation} that $\mathfrak{s}=\mathfrak{s}_{a,b}$.
		Conversely, for every $a\in\{0,\dots,\mathrm{ord}_p(n)\}$ and $b\in \mathbb{Z}$ the above argument also shows that
		$$\mathfrak{s}_{a,b}=\{r\in \mathfrak{R}~|~v(r-z)\ge d\}$$
		for $z=\zeta_n^b$ and $d=\frac{1}{p^{a-1}(p-1)}$. Hence, every $\mathfrak{s}_{a,b}$ indeed occurs as a cluster of $f_n$. Since $\{0\}$ and $\mathfrak{R}$ clearly occur as clusters of $f_n$, we have proved the first statement of the lemma.
		
		For the proof of the second statement, we first note that
		$$0\le d_{\mathfrak{R}}\le v(0-1)=v(1)=0$$
		since all elements in $\mathfrak{R}$ have non-negative valuations and $0,1\in\mathfrak{R}$.
		If $a\ge 1$, the cluster $\mathfrak{s}_{a,b}$ is proper for every $b\in\mathbb{Z}$. Using Lemma \ref{lem_cyclotomic-valuation}, we can compute the depth of $\mathfrak{s}_{a,b}$ by
		\begin{align*}
		d_{\mathfrak{s}_{a,b}}=\min_{k\neq k'}v\left(\zeta_n^{\frac{kn}{p^a}+b}-\zeta_n^{\frac{k'n}{p^a}+b}\right)=\min_{k\in\mathbb{Z}}v\left(\zeta_n^{\frac{kn}{p^a}}-1\right)=\frac{1}{p^{a-1}(p-1)}.
		\end{align*}
		If $a<\mathrm{ord}_p(n)$, then $\mathfrak{s}'=\mathfrak{s}_{a+1,b}$ is the smallest cluster satisfying $\mathfrak{s}_{a,b}\subsetneq \mathfrak{s}'$. Thus, we get
		$$\delta_{\mathfrak{s}_{a,b}}=d_{\mathfrak{s}_{a,b}}-d_{\mathfrak{s}_{a+1,b}}=\frac{1}{p^{a-1}(p-1)}-\frac{1}{p^{a}(p-1)}=\frac{1}{p^a}$$
		for the relative depth of $\mathfrak{s}_{a,b}$. If $a=\mathrm{ord}_p(n)$, then $\mathfrak{s}'=\mathfrak{R}$ is the smallest cluster satisfying $\mathfrak{s}_{a,b}\subsetneq \mathfrak{s}'$. Hence, we get
		$$\delta_{\mathfrak{s}_{a,b}}=d_{\mathfrak{s}_{a,b}}-d_{\mathfrak{R}}=\frac{1}{p^{a-1}(p-1)}$$
		for the relative depth of $\mathfrak{s}_{a,b}$. This completes the proof of the lemma.
	\end{proof}
	
	\section{Contributions of finite places}\label{sec_finite}
	We study the contributions of the finite places to the Faltings height. More precisely, we discuss the computation of $\mathrm{ord}_{\mathfrak{p}}(\Lambda)$ for any prime ideal $\mathfrak{p}$ and we give an explicit expression for the curve $X_n$.
	
	Let $\mathfrak{p}\subseteq \mathcal{O}_{K}$ be any non-zero prime ideal.
	If $\mathcal{X}$ has good reduction at $\mathfrak{p}$ then the base change $\Lambda_{\mathfrak{p}}$ of $\Lambda$ to $B_{\mathfrak{p}}$ is a trivializing section of the line bundle $\det H^0(\mathcal{X}_{\mathfrak{p}},\omega_{\mathcal{X}_{\mathfrak{p}}/B_\mathfrak{p}})^{\otimes (8g+4)}$, see \cite[Proposition 3.1]{dJo07}. Thus, we get
	\begin{align}\label{equ_good-reduction}
			\mathrm{ord}_{\mathfrak{p}}(\Lambda)=0
	\end{align}
	whenever $\mathcal{X}$ has good reduction at $\mathfrak{p}$. We can apply this to the hyperelliptic curve $X_n$ defined in Equation (\ref{equ_curve-Xn}) for an odd integer $n\ge 3$. We can choose the number field $K$ such that $X_n$ has semistable reduction over $B=\mathrm{Spec}\,\mathcal{O}_K$ and $\mathbb{Q}(i,\zeta_n,(-4)^{1/n})\subseteq K$. We write $\mathcal{X}_n$ for the minimal regular model and $\Lambda_n=\Lambda$ for the canonical rational section in $H^0(\mathcal{X}_n,\omega_{\mathcal{X}_n/B})^{\otimes (8g+4)}$ as above. We compute $\mathrm{ord}_{\mathfrak{p}}(\Lambda_n)$ in the following two lemmas.
	\begin{Lem}\label{lem_good-reduction}
		For any non-zero prime ideal $\mathfrak{p}\subseteq \mathcal{O}_K$ with $n\notin \mathfrak{p}$ it holds
		$$\mathrm{ord}_{\mathfrak{p}}(\Lambda_n)=0.$$
	\end{Lem}
	\begin{proof}
		By Equation (\ref{equ_good-reduction}) it is enough to show that $\mathcal{X}_n$ has good reduction at $\mathfrak{p}$. Since $\mathcal{X}_n$ is minimal, it is enough to find a Weierstrass equation for $X_n$, which defines a non-singular curve over $B_{\mathfrak{p}}$. By Remark \ref{rem_coordinate-change} the curve $X_n$ can be defined by the equation
		$$y^2+y=x^n$$
		The \emph{hyperelliptic discriminant} \cite[Definition 1.6]{Loc94} corresponding to this Weierstrass equation is given by
		\begin{align*}
			\Delta&=2^{2(n-1)}\mathrm{Disc}\left(x^{n}+\tfrac{1}{4}\right)=2^{2(n-1)}\left(\tfrac{1}{(-4)^{1/n}}\right)^{n(n-1)}\mathrm{Disc}\left(x^{n}-1\right)\\
			&=(-1)^{n(n-1)/2}n^{n},
		\end{align*}
		where we have used that the genus of $X_n$ is given by $g=(n-1)/2$.
		As $n\notin\mathfrak{p}$, we get $\Delta\not\equiv 0(\mathrm{mod}~\mathfrak{p})$.
		Thus, by \cite[Theorem 1.7]{Loc94}, the hyperelliptic curve defined by $y^2+y=x^{n}$ is non-singular over $\mathcal{O}_K/\mathfrak{p}$. It follows that $\mathcal{X}_n$ has good reduction at $\mathfrak{p}$.
	\end{proof}
	
	Next, we consider the case $n\in\mathfrak{p}$. Since $n$ is odd, it holds $2\notin \mathfrak{p}$. First, let $X$ be a hyperelliptic curve of genus $g>1$ defined by an equation $y^2=f(x)$ for a monic polynomial $f$ of degree $\deg f=2g+2$. We consider the differential
	$$\omega=\frac{dx}{2y}\wedge\frac{xdx}{2y}\wedge\dots\wedge\frac{x^{g-1}dx}{2y}.$$
	Kunzweiler \cite[Theorem 1.4]{Kun20} showed that if $\lambda\cdot \omega$ is a basis of $\det H^0(\mathcal{X}_\mathfrak{p},\omega_{\mathcal{X}_\mathfrak{p}/B_{\mathfrak{p}}})$, then
	$$8v(\lambda)=\sum_{\substack{|\mathfrak{s}|\text{ even}\\ \mathfrak{s}\neq \mathfrak{R}}}\delta_{\mathfrak{s}}(|\mathfrak{s}|-2)|\mathfrak{s}|+\sum_{|\mathfrak{s}|\text{ odd}}\delta_{\mathfrak{s}}(|\mathfrak{s}|-1)^2+d_{\mathfrak{R}}2g(2g+2).$$
	Here, the sums run over all proper clusters corresponding to the polynomial $f$, where the first sum only takes clusters of even cardinality and unequal to $\mathfrak{R}$ into account and the second sum only takes clusters of odd cardinality into account.
	Since $\Lambda=2^{4g^2}\mathrm{Disc}(f)^g\omega^{\otimes (8g+4)}$ and $2\notin \mathfrak{p}$, we can conclude
	\begin{align}\label{equ_integral-differential}
		\frac{\mathrm{ord}_{\mathfrak{p}}(\Lambda)}{e_\mathfrak{p}}=&g\cdot v(\mathrm{Disc}(f))-\tfrac{2g+1}{2}\sum_{\substack{|\mathfrak{s}|\text{ even}\\ \mathfrak{s}\neq \mathfrak{R}}}\delta_{\mathfrak{s}}(|\mathfrak{s}|-2)|\mathfrak{s}|\\
		&-\tfrac{2g+1}{2}\sum_{|\mathfrak{s}|\text{ odd}}\delta_{\mathfrak{s}}(|\mathfrak{s}|-1)^2-d_{\mathfrak{R}}g(2g+2)(2g+1).\nonumber
	\end{align}
	
	\begin{Lem}\label{lem_ord-reduction}
		For any odd integer $n\ge 5$ and any non-zero prime ideal $\mathfrak{p}\subseteq \mathcal{O}_K$ with $n\in \mathfrak{p}$ it holds
		$$\mathrm{ord}_{\mathfrak{p}}(\Lambda_n)=\frac{n}{2}\left(n\frac{p^{2\mathrm{ord}_{p}(n)}-1}{p^{2\mathrm{ord}_{p}(n)-1}(p^2-1)}-\mathrm{ord}_{p}(n)\right)e_{\mathfrak{p}},$$
		where $p$ denotes the unique prime number with $p\in\mathfrak{p} \cap\mathbb{Z}$.
	\end{Lem}
	\begin{proof}
		The proof directly follows by applying Lemma \ref{lem_cluster} to Equation (\ref{equ_integral-differential}). For an explicit calculation, let us briefly write $n_p=\mathrm{ord}_p(n)$. Then we get
		$$v(\mathrm{Disc}(f_n))=v(\mathrm{Disc}(x(1-x^{n})))=v(n^n)=n\cdot n_p$$
		for the valuation of the discriminant of the polynomial $f_n(x)=x(1-x^{n})$, 
		$$\sum_{\substack{|\mathfrak{s}|\text{ even}\\ \mathfrak{s}\neq \mathfrak{R}}}\delta_{\mathfrak{s}}(|\mathfrak{s}|-2)|\mathfrak{s}|=0$$
		as by Lemma \ref{lem_cluster} the cardinality $|\mathfrak{s}_{a,b}|=p^a$ is always odd, and
		\begin{align*}
		&\sum_{|\mathfrak{s}|\text{ odd}}\delta_{\mathfrak{s}}(|\mathfrak{s}|-1)^2=\sum_{a=1}^{n_p-1}\sum_{b=1}^{n/p^a}\frac{(p^a-1)^2}{p^a}+\sum_{b=1}^{n/p^{n_p}}\frac{(p^{n_p}-1)^2}{p^{n_p-1}(p-1)}\\
		&=n\left(\sum_{a=1}^{n_p-1}\left(1-\frac{2}{p^a}+\frac{1}{p^{2a}}\right)+\frac{1}{p-1}\left(p-\frac{2}{p^{n_p-1}}+\frac{1}{p^{2n_p-1}}\right)\right)\\
		&=n\left(n_p-\frac{p^{2n_p}-1}{p^{2n_p-1}(p^2-1)}\right).
		\end{align*}
		Since $d_\mathfrak{R}=0$ by Lemma \ref{lem_cluster} and $g=\frac{n-1}{2}$, we can conclude by Equation (\ref{equ_integral-differential})
		\begin{align*}
			\frac{\mathrm{ord}_{\mathfrak{p}}(\Lambda_n)}{e_{\mathfrak{p}}}&=\frac{n-1}{2}\cdot n\cdot n_p-\frac{n}{2}\cdot n\cdot \left(n_p-\frac{p^{2n_p}-1}{p^{2n_p-1}(p^2-1)}\right)\\
			&=\frac{n}{2}\left(n\cdot \frac{p^{2n_p}-1}{p^{2n_p-1}(p^2-1)}-n_p\right).
		\end{align*}
		This proves the lemma.
	\end{proof}

	\section{Hyperelliptic period matrices}\label{sec_period-matrix}
	To compute the contributions of the infinite places to the Faltings height in the next section, we study the period matrix of a hyperelliptic Riemann surface in this section.
	
	Let $M$ be any compact and connected Riemann surface of genus $g\ge 1$, which is hyperelliptic. We can present $M$ as the closure in a weighted projective space of an equation
	$$y^2=f(x):=\prod_{i=1}^{m}(x-a_i)$$
	for pairwise distinct $a_1,\dots,a_{m}\in\mathbb{C}$, where $m=2g+1$ or $m=2g+2$. For $1\le j\le g$ we write $\omega_j=\frac{x^{j-1}dx}{y}$. Then $\omega_1,\dots,\omega_g$ form a basis of $H^0(M,\Omega_M^1)$, see for example \cite[Proposition 1.12]{Loc94}.
	As in \cite[Chapter IIIa, §5]{Mum84}, we choose a symplectic basis for homology $H_1(M,\mathbb{Z})$ of $M$ in the following way: On $\mathbb{P}_{\mathbb{C}}^1$ we choose curves $A_i$ circling the points $a_{2i-1}$ and $a_{2i}$ and curves $B_i$ circling the points $a_{2i},a_{2i+1},\dots,a_{2g+1}$ for all $1\le i\le g$. We can choose these curves, such that the $A_i$'s are pairwise disjoint, the $B_i$'s are pairwise disjoint and $A_i$ and $B_j$ only intersect if $i=j$. The Riemann surface $M$ admits a double cover $M\to\mathbb{P}_{\mathbb{C}}^1$ sending $(x,y)$ to $x$, which ramifies exactly at the points $a_1,\dots,a_{2g+2}$, where $a_{2g+2}=\infty$ if $\deg f=2g+1$. As the curves $A_i$ and $B_i$ circle an even number of ramification points, they can be lifted to $M$. We denote their lifts also by $A_i$ and $B_i$. They form a symplectic basis for homology $H_1(M,\mathbb{Z})$ of $M$.
	The corresponding period matrix $\Omega=(\Omega_1|\Omega_2)\in\mathbb{C}^{g\times 2g}$ is defined by $(\Omega_1)_{ij}=\int_{A_j}\omega_i$ and $(\Omega_2)_{ij}=\int_{B_j}\omega_i$.
	
	Let us make the entries of $\Omega$ more explicit. The double cover $M\to\mathbb{P}_{\mathbb{C}}^1$ presents $M$ as two layers of $\mathbb{P}_{\mathbb{C}}^1$. We can think of the integral over the cycle $A_k$ as the sum of the integrals from $a_{2k-1}$ to $a_{2k}$ on the first layer and from $a_{2k}$ to $a_{2k-1}$ on the second layer. This means
	$$(\Omega_1)_{jk}=\int_{A_k}\frac{x^{j-1}dx}{y}=\int_{a_{2k-1}}^{a_{2k}}\frac{x^{j-1}dx}{\sqrt{f(x)}}+\int_{a_{2k}}^{a_{2k-1}}\frac{x^{j-1}dx}{-\sqrt{f(x)}}=2\int_{a_{2k-1}}^{a_{2k}}\frac{x^{j-1}dx}{\sqrt{f(x)}}.$$
	Similarly, the integrals over the cycles $B_k$ are given by
	$$(\Omega_2)_{jk}=2\int_{a_{2k}}^{a_{2g+1}}\frac{x^{j-1}dx}{\sqrt{f(x)}}.$$
	
	As shown in \cite[p.~231]{GH78}, we have
	\begin{align}\label{equ_pairing-matrix}
	\int_M\omega\wedge\overline{\omega'}=\sum_{j=1}^g\left(\int_{A_j}\omega\cdot\overline{\int_{B_j}\omega'}-\int_{B_j}\omega\cdot\overline{\int_{A_j}\omega'}\right)
	\end{align}
	for any holomorphic sections $\omega,\omega'\in H^0(M,\Omega_M^1)$. Since the pairing $\frac{i}{2}\int_M\omega\wedge\overline{\omega'}$ is positive-definite on $H^0(M,\Omega_M^1)$, we in particular obtain that the linear form 
	$$\mathbb{C}^g\to \mathbb{C}^g,\qquad v\mapsto\ltrans{\Omega_2} v=\left(\int_{B_j}\sum_{i=1}^gv_i\omega_i\right)_{1\le j\le g}$$ has a trivial kernel and thus $\Omega_2$ is invertible. Applying Equation (\ref{equ_pairing-matrix}) to the basis $\omega_1,\dots,\omega_g$, we deduce that
	$$\left(\int_M\omega_j\wedge\overline{\omega_k}\right)_{1\le j,k\le g}=\Omega_1\cdot\ltrans{\overline{\Omega}_2}-\Omega_2\cdot\ltrans{\overline{\Omega}_1}.$$
	Note that by the first Riemann bilinear relation \cite[p.~231]{GH78} we have $\Omega_1\cdot\ltrans{\Omega}_2=\Omega_2\cdot\ltrans{\Omega}_1$.
	Thus, we can compute
	\begin{align*}
		&\det\left(\int_M\omega_j\wedge\overline{\omega_k}\right)_{1\le j,k\le g}=\det(\Omega_1\cdot\ltrans{\overline{\Omega}_2}-\Omega_2\cdot\ltrans{\overline{\Omega}_1})\\
		&=\det(\Omega_1\cdot\ltrans{\overline{\Omega}_2}-\Omega_2\cdot\overline{\Omega}_2^{-1}\cdot \overline{\Omega}_2\cdot \ltrans{\overline{\Omega}_1})=\det(\Omega_1\cdot\ltrans{\overline{\Omega}_2}-\Omega_2\cdot\overline{\Omega}_2^{-1}\cdot \overline{\Omega}_1\cdot \ltrans{\overline{\Omega}_2})\\
		&=\det(\Omega_1-\Omega_2\cdot\overline{\Omega}_2^{-1}\cdot\overline{\Omega}_1)\cdot \det(\overline{\Omega}_2)=\det \begin{pmatrix} \Omega_1&\Omega_2\\ \overline{\Omega}_1& \overline{\Omega}_2\end{pmatrix}.
	\end{align*}
	Using linear dependencies, we get
	\begin{align}\label{equ_determinant-A}
		\det\left(\int_M\omega_j\wedge\overline{\omega_k}\right)_{1\le j,k\le g}=\det\begin{pmatrix} \Omega_1&\Omega_2\\ \overline{\Omega}_1& \overline{\Omega}_2\end{pmatrix}=\det\begin{pmatrix} A\\ \overline{A}\end{pmatrix}
	\end{align}
	for the $g\times 2g$ matrix $A$ with entries
	$$A_{jk}=2\int_{a_{2g+1}}^{a_k}\frac{x^{j-1}dx}{\sqrt{f(x)}}.$$
	After a coordinate change, we may assume that $a_{2g+1}=0$. Then we get
	\begin{align}\label{equ_matrix-entries}
		A_{jk}=2\int_0^{a_k}\frac{x^{j-1}dx}{\sqrt{f(x)}}=2a_k^{j}\int_0^1\frac{x^{j-1}dx}{\sqrt{f(a_k x)}}.
	\end{align}

	\section{Contributions of infinite places}\label{sec_infinite}
	We study the contributions of the infinite places to the Faltings height. If $X$ is a hyperelliptic curve defined over a number field $K$, we get a Riemann surface $X_{\sigma}$ by the base change along any complex embedding $\sigma\colon K\to\mathbb{C}$. We will explicitly compute the Faltings metric of the canonical element $\Lambda$ of $\det H^0(M,\Omega_M^1)^{\otimes (8g+4)}$ for the hyperelliptic Riemann surface $M$ given by the equation $y^2=x(1-x^n)$ for any $n\ge 2$.
	
	For a hyperelliptic Riemann surface $M$ as in the previous section, the Faltings metric on $H^0(M,\Omega_M^1)$ is given by the inner product
	$$\langle\eta,\eta'\rangle=\frac{i}{2}\int_M \eta\wedge \overline{\eta'}$$
	for any $\eta,\eta'\in H^0(M,\Omega_M^1)$. It induces the inner product
	$$\langle \eta_1\wedge\dots\wedge\eta_g,\eta_1'\wedge\dots\wedge\eta_g'\rangle=\det\left(\frac{i}{2}\int_M\eta_j\wedge\overline{\eta_k'}\right)_{1\le j,k\le g}$$
	on $\det H^0(M,\Omega_M^1)$. The canonical section $\Lambda\in \det H^0(M,\Omega_M^1)^{\otimes (8g+4)}$ has the form 
	\begin{align}\label{equ_complex-lambda}
		\Lambda=\frac{\mathrm{Disc}(f)^g}{2^{4g(g+1)}}(\omega_1\wedge\dots\wedge\omega_g)^{\otimes (8g+4)}.
	\end{align}
	It follows that its Faltings metric is given by
	\begin{align}\label{equ_lambda-norm}
	\|\Lambda\|^2=\frac{|\mathrm{Disc}(f)|^{2g}}{2^{4g(4g+3)}}\left(\det\left(\int_M\omega_j\wedge\overline{\omega_k}\right)_{1\le j,k\le g}\right)^{8g+4}.
	\end{align}
	
	Now we consider the hyperelliptic Riemann surface $M_n$ associated to the equation
	$$y^2=x(1-x^n)$$
	for any $n\ge 2$. We denote $\Lambda_n\in\det H^0(M_n,\Omega_{M_n}^1)^{\otimes (8g+4)}$ for the corresponding canonical element in Equation (\ref{equ_complex-lambda}). The following lemma computes its Faltings metric. We have also included the case of even $n$ for completeness.
	\begin{Lem}\label{lem_faltings-metric}
		We have $\|\Lambda_{n}\|=\frac{\pi^{g(4g+2)}}{2^{2g}n^{n(g+1)}}\left(\prod_{j=1}^g\frac{\Gamma\left(\frac{2j-1}{2n}\right)}{\Gamma\left(\frac{1}{2}+\frac{2j-1}{2n}\right)}\right)^{8g+4}$ for every integer $n\ge 2$, where $g=\left\lfloor \frac{n}{2}\right\rfloor$ denotes the genus of $M_n$.
	\end{Lem}
	\begin{proof}
		We have $n=2g+1$ if $n$ is odd and $n=2g$ if $n$ is even. We write $\zeta_n=e^{2\pi i/n}$ for the $n$-th root of unity. We can set $a_j=\zeta_n^j$ for $1\le j\le 2g$ and $a_{2g+1}=0$. If $n$ is odd, we get $a_{2g+2}=1$. If $n$ is even, we get $a_{2g+2}=\infty$.
		We can compute the value of (\ref{equ_matrix-entries}) in this case explicitly by
		$$A_{jk}=2\zeta_n^{k\left(j-\frac{1}{2}\right)}\int_0^1x^{j-\frac{3}{2}}(1-x^n)^{-\frac{1}{2}}dx=\frac{2\zeta_n^{k\left(j-\frac{1}{2}\right)}}{n}\int_0^1 t^{\frac{2j-1}{2n}-1}(1-t)^{-\frac{1}{2}}dt.$$
		Using Euler's beta function, we can compute the integral as
		$$\int_0^1 t^{\frac{2j-1}{2n}-1}(1-t)^{-\frac{1}{2}}dt=B\left(\tfrac{2j-1}{2n},\tfrac{1}{2}\right)=\frac{\Gamma\left(\frac{2j-1}{2n}\right)\Gamma\left(\frac{1}{2}\right)}{\Gamma\left(\frac{1}{2}+\frac{2j-1}{2n}\right)}.$$
		Note that $\Gamma(\frac{1}{2})=\sqrt{\pi}$. Now we can compute the determinant in Equation (\ref{equ_determinant-A}) by
		$$\det\begin{pmatrix}A\\ \overline{A}\end{pmatrix}= \left(\frac{2^{g}\pi^{g/2}\prod_{j=1}^{g}\Gamma\left(\frac{2j-1}{2n}\right)}{n^{g}\prod_{j=1}^{g}\Gamma\left(\frac{n+2j-1}{2n}\right)}\right)^2\det\begin{pmatrix}\left(\zeta_n^{k\left(j-\frac{1}{2}\right)}\right)_{1\le j\le g,1\le k\le 2g} \\ \left(\overline{\zeta}_n^{k\left(j-\frac{1}{2}\right)}\right)_{1\le j\le g,1\le k\le 2g}\end{pmatrix}.$$
		The determinant in the last expression can be computed via a Vandermonde matrix
		\begin{align*}
			&\det\begin{pmatrix}\left(\zeta_n^{k\left(j-\frac{1}{2}\right)}\right)_{1\le j\le g,1\le k\le 2g} \\ 	\left(\overline{\zeta}_n^{k\left(j-\frac{1}{2}\right)}\right)_{1\le j\le g,1\le k\le 2g}\end{pmatrix}=\det\begin{pmatrix}\left(\zeta_n^{k\left(j-\frac{1}{2}\right)}\right)_{1\le j\le g,0\le k\le 2g-1} \\ \left(\zeta_n^{k\left(n-j+\frac{1}{2}\right)}\right)_{1\le j\le g,0\le k\le 2g-1}\end{pmatrix}\\
			&=\zeta_n^{-\frac{g(2g-1)}{2}}\det\begin{pmatrix}\left(\zeta_n^{kj}\right)_{1\le j\le g,0\le k\le 2g-1} \\ \left(\zeta_n^{k\left(n-j+1\right)}\right)_{1\le j\le g,0\le k\le 2g-1}\end{pmatrix}\\
			&=\zeta_n^{-\frac{g(2g-1)}{2}}(-1)^{\frac{g(g-1)}{2}}\det\begin{pmatrix}\left(\zeta_n^{kj}\right)_{1\le j\le g,0\le k\le 2g-1} \\ \left(\zeta_n^{k\left(j+n-g\right)}\right)_{1\le j\le g,0\le k\le 2g-1}\end{pmatrix}\\
			&=\begin{cases}\zeta_n^{-\frac{g(2g-1)}{2}}(-1)^{\frac{g(g-1)}{2}}\prod_{1\le j<k\le n}(\zeta_n^k-\zeta_n^j) & \text{if } n=2g,\\ \zeta_n^{-\frac{g(2g-1)}{2}}(-1)^{\frac{g(g-1)}{2}}\frac{\prod_{1\le j<k\le n}(\zeta_n^k-\zeta_n^j)}{\prod_{j=1}^g(\zeta_n^{g+1}-\zeta_n^j)\prod_{j=g+2}^{n}(\zeta_n^j-\zeta_n^{g+1})} & \text{if } n=2g+1. 
			\end{cases}
		\end{align*}
		Note that in any case we have
		$$\prod_{1\le j<k\le n}(\zeta_n^k-\zeta_n^j)^2=|\mathrm{Disc}(x^n-1)|=n^n.$$
		For the denominator in the case $n=2g+1$ we get
		$$\prod_{j=1}^g(\zeta_n^{g+1}-\zeta_n^j)\prod_{j=g+2}^{n}(\zeta_n^j-\zeta_n^{g+1})=\zeta_n^{2g(g+1)}(-1)^g\prod_{j=1}^{2g}(1-\zeta_n^{j})=n\zeta_n^{g}(-1)^g,$$
		where the last equality follows since
		$$\prod_{j=1}^{2g}(X-\zeta_n^{j})=\frac{X^n-1}{X-1}=\sum_{j=0}^{2g}X^j.$$
		Thus, taking the $(8g+4)$-th power, all signs cancel out. We get by putting everything together and using Equation (\ref{equ_determinant-A})
		$$	\det\left(\int_M\omega_j\wedge\overline{\omega_k}\right)_{1\le j,k\le g}^{8g+4}=n^{-n(4g+2)} \left(\frac{2^{g}\pi^{g/2}\prod_{j=1}^{g}\Gamma\left(\frac{2j-1}{2n}\right)}{\prod_{j=1}^{g}\Gamma\left(\frac{1}{2}+\frac{2j-1}{2n}\right)}\right)^{16g+8}$$
		both if $n$ is even and if $n$ is odd. Applying this to Equation (\ref{equ_lambda-norm}) we obtain
		$$\|\Lambda_{n}\|=\frac{\pi^{g(4g+2)}}{2^{2g}n^{n(g+1)}}\left(\prod_{j=1}^g\frac{\Gamma\left(\frac{2j-1}{2n}\right)}{\Gamma\left(\frac{1}{2}+\frac{2j-1}{2n}\right)}\right)^{8g+4}.$$
		This completes the proof of the lemma.
	\end{proof}
	
	\section{Proof of Theorem \ref{thm_faltings-height-hyperelliptic}}\label{sec_proof-main-thm}
	We can now put the computations from the previous sections together to prove Theorem \ref{thm_faltings-height-hyperelliptic}. First, we remark that the case $n=3$ was already shown by Deligne \cite[p. 29]{Del85}. Indeed, Deligne showed the formula
	$$h'_{\mathrm{Fal}}(X_3)=-\tfrac{1}{2}\log\frac{\Gamma(1/3)^3}{\sqrt{3}\cdot\Gamma(2/3)^3}.$$
	Note that he used the different normalization $h'_{\mathrm{Fal}}(X)=h_{\mathrm{Fal}}(X)+\frac{g}{2}\log\pi$ of the Faltings height. Using the classical formulas for the gamma function
	$$\Gamma(1/6)=\frac{\sqrt{\frac{3}{\pi}}\cdot\Gamma(1/3)^2}{\sqrt[3]{2}}\qquad \text{and}\qquad \Gamma(1/3)\cdot\Gamma(2/3)=\frac{2\pi}{\sqrt{3}},$$
	one checks that this coincides with the formula in Theorem \ref{thm_faltings-height-hyperelliptic}.
	
	Thus, we can assume in the following that $n$ is an odd integer satisfying $n\ge 5$. Similarly to Equation (\ref{equ_heightbysection}), we will use the canonical section $\Lambda_n$ of the determinant line bundle $\det H^0(\mathcal{X}_n,\omega_{\mathcal{X}_n/B})^{\otimes (8g+4)}$ to compute the Faltings height by
	\begin{align*}
		(8g+4)h_{\mathrm{Fal}}(X_n)&=\tfrac{1}{d}\widehat{\deg}\det H^0(\mathcal{X}_n,\omega_{\mathcal{X}_n/B})^{\otimes(8g+4)}\\
		&=\tfrac{1}{d}\sum_{\mathfrak{p}\in|B|}\mathrm{ord}_{\mathfrak{p}}(\Lambda_n)\log N(\mathfrak{p})-\tfrac{1}{d}\sum_{\sigma\colon K\to\mathbb{C}}\log\|\Lambda_n\|_\sigma.	
	\end{align*}
	where $d=[K:\mathbb{Q}]$.
	We apply Lemmas \ref{lem_good-reduction} and \ref{lem_ord-reduction} to the vanishing orders in the first sum and Lemma \ref{lem_faltings-metric} to the norm in the second sum. This yields
	\begin{align*}
		(8g+4)h_{\mathrm{Fal}}(X_n)=&\tfrac{1}{d}\sum_{\substack{p\text{ prime}\\ p\mid n}}\sum_{\substack{\mathfrak{p}\in |B|\\ \mathfrak{p}\mid p}} \frac{n}{2}\left(n\frac{p^{2\mathrm{ord}_{p}(n)}-1}{p^{2\mathrm{ord}_{p}(n)-1}(p^2-1)}-\mathrm{ord}_{p}(n)\right)e_{\mathfrak{p}}\log N(\mathfrak{p})\\
		&-\tfrac{1}{d}\sum_{\sigma\colon K\to\mathbb{C}}\log\left(\frac{\pi^{g(4g+2)}}{2^{2g}n^{n(g+1)}}\left(\prod_{j=1}^g\frac{\Gamma\left(\frac{2j-1}{2n}\right)}{\Gamma\left(\frac{1}{2}+\frac{2j-1}{2n}\right)}\right)^{8g+4}\right).
	\end{align*}
	Since $\prod_{\mathfrak{p}\mid p}N(\mathfrak{p})^{e_{\mathfrak{p}}}=N\left(\prod_{\mathfrak{p}\mid p}\mathfrak{p}^{e_{\mathfrak{p}}}\right)=N(p\mathcal{O}_K)=p^d$ and since there are exactly $d$ embeddings $\sigma\colon K\to\mathbb{C}$, we can simplify the term to
	\begin{align*}
	(8g+4)h_{\mathrm{Fal}}(X_n)=&\sum_{\substack{p\text{ prime}\\ p\mid n}}\frac{n}{2}\left(n\frac{p^{2\mathrm{ord}_{p}(n)}-1}{p^{2\mathrm{ord}_{p}(n)-1}(p^2-1)}-\mathrm{ord}_{p}(n)\right)\log p\\
	&-\log \frac{\pi^{g(4g+2)}}{2^{2g}n^{n(g+1)}}-(8g+4)\sum_{j=1}^g\log\frac{\Gamma\left(\frac{2j-1}{2n}\right)}{\Gamma\left(\frac{1}{2}+\frac{2j-1}{2n}\right)}.
	\end{align*}
	Dividing both sides by $4n=8g+4$ and reordering the terms, we get
	\begin{align*}
		h_{\mathrm{Fal}}(X_n)=&\tfrac{n}{8}\sum_{\substack{p\text{ prime}\\ p\mid n}}\frac{p^{2\mathrm{ord}_{p}(n)}-1}{p^{2\mathrm{ord}_{p}(n)-1}(p^2-1)}\log p+\tfrac{n}{8}\log n-\tfrac{g}{2}\log \pi+\tfrac{g}{2n}\log 2\\
		&-\sum_{j=1}^g\log\frac{\Gamma\left(\frac{2j-1}{2n}\right)}{\Gamma\left(\frac{1}{2}+\frac{2j-1}{2n}\right)}.
	\end{align*}
	This completes the proof of Theorem \ref{thm_faltings-height-hyperelliptic}.
	
	\section{Bounds for products of gamma values}\label{sec_bound}
	We give bounds for the product of values of the gamma function occurring in Theorem \ref{thm_faltings-height-hyperelliptic}. More precisely, we will compare the sum of the logarithms of the values of the gamma function in Theorem \ref{thm_faltings-height-hyperelliptic} with the integral of the logarithm of the gamma function. For this purpose, we recall that
	\begin{align}\label{equ_integration-of-gamma}
	\int_0^z\log \Gamma(x)dx=\frac{z(1-z)}{2}+\frac{z}{2}\log 2\pi +z\log \Gamma(z)-\log G(1+z),
	\end{align}
	where $G$ denotes the Barnes $G$-function. We recall that the Barnes $G$-function satisfies the functional equation $G(z+1)=G(z)\Gamma(z)$ with the normalization $G(1)=1$.
	The goal of this section is to prove the following lemma.
	\begin{Lem}\label{lem_bound-gamma}
		For every odd integer $n\ge 3$, we have the lower bound
		$$\sum_{j=1}^{g}\log\frac{\Gamma\left(\frac{2j-1}{2n}\right)}{\Gamma\left(\frac{1}{2}+\frac{2j-1}{2n}\right)}\ge \frac{3-\log 2-36\zeta'(-1)}{12}n-\frac{1+\log 2n}{2}$$
		and the upper bound
		$$\sum_{j=1}^{g}\log\frac{\Gamma\left(\frac{2j-1}{2n}\right)}{\Gamma\left(\frac{1}{2}+\frac{2j-1}{2n}\right)}\le\frac{3-\log 2-36\zeta'(-1)}{12}n,$$
		where $g=(n-1)/2$.
	\end{Lem}
	\begin{proof}
		We briefly write
		$$I_a^b=\int_a^b\log\Gamma(x)dx$$
		for any $b\ge a\ge 0$.	For any $z\in\left[\frac{1}{2n},1-\frac{1}{n}\right]$ we have
		$$nI_z^{z+\frac{1}{n}}\le \log\Gamma(z)\le nI_{z-\frac{1}{2n}}^{z+\frac{1}{2n}}.$$
		The first inequality follows as $\log \Gamma(x)$ is monotonic decreasing on $(0,1]$. The second inequality follows by Jensen's inequality as $\log \Gamma(x)$ is convex on $(0,\infty)$. Thus, by summing up we obtain
		$$nI_{\frac{1}{2n}}^{\frac{1}{2}}-nI_{\frac{1}{2}}^{1-\frac{1}{2n}}\le \sum_{j=1}^g\log\frac{\Gamma\left(\frac{2j-1}{2n}\right)}{\Gamma\left(\frac{1}{2}+\frac{2j-1}{2n}\right)}\le nI_0^{\frac{1}{2}-\frac{1}{2n}}-nI_{\frac{1}{2}+\frac{1}{2n}}^{1}.$$
		Using the non-negativity and the monotonicity of $\log\Gamma$ on $(0,1]$, we can further estimate
		$$n\left(2I_0^{\frac{1}{2}}-I_0^{1}-I_0^{\frac{1}{2n}}\right)\le \sum_{j=1}^g\log\frac{\Gamma\left(\frac{2j-1}{2n}\right)}{\Gamma\left(\frac{1}{2}+\frac{2j-1}{2n}\right)}\le n\left(2I_0^{\frac{1}{2}}-I_0^1\right).$$
		By Equation (\ref{equ_integration-of-gamma}) the values $I_0^1$ and $I_0^{\frac{1}{2}}$ are given by
		$$I_0^1=\tfrac{1}{2}\log 2\pi,\qquad I_0^{\frac{1}{2}}=\tfrac{1}{8}+\tfrac{1}{4}\log 2-\log G\left(\tfrac{1}{2}\right),$$
		where we used $\Gamma(1)=G(2)=1$ and $\Gamma(1/2)=\sqrt{\pi}$.
		Since $G\left(\frac{1}{2}\right)=2^{\frac{1}{24}}e^{\frac{3}{2}\zeta'(-1)}\pi^{-\frac{1}{4}}$, we get
		$$2I_0^{\frac{1}{2}}-I_0^1=\tfrac{1}{4}-\tfrac{1}{12}\log 2-3\zeta'(-1).$$
		Since $\Gamma(x)\cdot x=\Gamma(x+1)$ and $\Gamma(t)\le 1$ for all $t\in[1,2]$, we get $\log\Gamma(x)\le -\log x$ for all $x\in(0,1]$. Hence, we get for $I_0^{\frac{1}{2n}}$ the upper bound
		$$I_0^{\frac{1}{2n}}\le \int_0^{\frac{1}{2n}}(-\log x)dx=\frac{1+\log 2n}{2n}.$$
		Putting everything together, we obtain the bounds claimed in the lemma.
	\end{proof}
	
	\section{Proof of Corollary \ref{cor_bounds}}\label{sec_proof-corollary}
	We bound the terms in the formula for the Faltings height in Theorem \ref{thm_faltings-height-hyperelliptic} to prove Corollary \ref{cor_bounds}. After establishing the bound for the values of the gamma function in the previous section, it remains to find a good bound for the sum over the primes dividing $n$ in the expression for the Faltings height of $X_n$. This is done by the following lemma.
	\begin{Lem}\label{lem_bound-finite}
		For every odd integer $n\ge 3$ it holds
		$$\frac{\log n}{n}\le\sum_{\substack{p\text{ prime}\\ p\mid n}}\frac{p^{2\mathrm{ord}_p(n)}-1}{p^{2\mathrm{ord}_p(n)-1}(p^2-1)}\log p< \tfrac{9}{8}\log\log n+0.7405.$$
	\end{Lem}
	\begin{proof}
		We first show the lower bound. We have
		$$\sum_{\substack{p\text{ prime}\\ p\mid n}}\frac{p^{2\mathrm{ord}_p(n)}-1}{p^{2\mathrm{ord}_p(n)-1}(p^2-1)}\log p\ge\sum_{\substack{p\text{ prime}\\ p\mid n}}\frac{\log p}{p}.$$
		Since the function $\frac{\log x}{x}$ is monotonically decreasing for $x\ge e$ and since $p\ge3$, we obtain
		$$\sum_{\substack{p\text{ prime}\\ p\mid n}} \frac{\log p}{p}\ge\sum_{\substack{p\text{ prime}\\ p\mid n}}\frac{\log p^{\mathrm{ord}_p(n)}}{p^{\mathrm{ord}_p(n)}}\ge \sum_{\substack{p\text{ prime}\\ p\mid n}}\frac{\log p^{\mathrm{ord}_p(n)}}{n}=\frac{\log n}{n}.$$
		This shows the lower bound of the lemma.
		
		Next, we prove the upper bound. Since $n\ge3$ is an odd integer, every prime factor $p$ of $n$ is at least $3$ and hence,
		\begin{align}\label{equ_bound-proof-primes}
		\sum_{\substack{p\text{ prime}\\ p\mid n}}\frac{p^{2\mathrm{ord}_p(n)}-1}{p^{2\mathrm{ord}_p(n)-1}(p^2-1)}\log p= \sum_{\substack{p\text{ prime}\\ p\mid n}}\frac{1-p^{-2\mathrm{ord}_p(n)}}{p\left(1-\frac{1}{p^2}\right)}\log p
		\le \tfrac{9}{8}\sum_{\substack{p\text{ prime}\\ p\mid n}}\frac{\log p}{p},
		\end{align}		
		We denote $\omega(n)$ for the number of distinct prime factors of $n$ and we write $p_i$ for the $i$-th prime number, that is $p_1=2$, $p_2=3$, $p_3=5$, and so on. Since $\frac{\log x}{x}$ is monotonically decreasing for $x\ge e$, we can estimate
		$$\sum_{\substack{p\text{ prime}\\ p\mid n}}\frac{\log p}{p}\le \sum_{i=2}^{\omega(n)+1}\frac{\log p_i}{p_i}.$$
		Using a result by Rosser and Schoenfeld \cite[Equation (3.24)]{RS62}, we can bound
		$$\sum_{i=2}^{\omega(n)+1}\frac{\log p_i}{p_i}=\sum_{\substack{p\text{ prime}\\ p\le p_{\omega(n)+1}}}\frac{\log p}{p}-\frac{\log 2}{2}\le \log p_{\omega(n)+1} -\frac{\log 2}{2}.$$
		Since $n$ is odd, we have $\omega(n)+1=\omega(2n)$. Using the Chebyshev function $\vartheta(x)$, we get
		$$0.5972p_{\omega(2n)}\le \vartheta(p_{\omega(2n)})=\sum_{i=1}^{\omega(2n)}\log p_i\le \log 2n,$$
		where the first inequality follows from \cite[Theorem 10]{RS62} for $p_{\omega(2n)}\ge 101$ and by a direct computation for $3\le p_{\omega(2n)}<101$. Since $n\ge 3$, we have $\log 2n\le 1.631\log n$.
		Thus, we conclude
		$$\sum_{\substack{p\text{ prime}\\ p\mid n}}\frac{\log p}{p}\le \log p_{\omega(2n)}-\frac{\log 2}{2}\le \log \frac{1.631\log n}{0.5972}-\frac{\log 2}{2}<\log\log n+0.6582.$$
		If we apply this to the bound in (\ref{equ_bound-proof-primes}), we obtain
		$$\sum_{\substack{p\text{ prime}\\ p\mid n}}\frac{p^{2\mathrm{ord}_p(n)}-1}{p^{2\mathrm{ord}_p(n)-1}(p^2-1)}\log p\le \tfrac{9}{8}\sum_{\substack{p\text{ prime}\\ p\mid n}}\frac{\log p}{p}< \tfrac{9}{8}\log\log n+0.7405.$$
		This completes the proof of the lemma.
	\end{proof}
	
	Now we can give the proof of Corollary \ref{cor_bounds}.
	\begin{proof}[Proof of Corollary \ref{cor_bounds}]
		If we apply Lemmas \ref{lem_bound-gamma} and \ref{lem_bound-finite} to the expression of $h_{\mathrm{Fal}}(X_n)$ in Theorem \ref{thm_faltings-height-hyperelliptic}, we obtain the lower bound
		\begin{align*}
		h_{\mathrm{Fal}}(X_n)&\ge \tfrac{n}{8}\log n-\tfrac{3-\log 2-36\zeta'(-1)}{12}n-\tfrac{g}{2}\log \pi+\tfrac{1}{8}\log n+\tfrac{g}{2n}\log 2\\
		&\ge \tfrac{n}{8}\log n-\tfrac{3-\log 2-36\zeta'(-1)+3\log \pi}{12}n\\
		&> \tfrac{n}{8}\log n-0.975n
		\end{align*}
		and the upper bound
		\begin{align}\label{equ_upper-bound}
		h_{\mathrm{Fal}}(X_n)\le &\tfrac{n}{8}\log n +\tfrac{9}{64}n\log\log n+\tfrac{0.7405}{8}n-\tfrac{3-\log2-36\zeta'(-1)}{12}n\\
		&-\tfrac{g}{2}\log \pi+\tfrac{\log 2n}{2}+\tfrac{g}{2n}\log 2+\tfrac{1}{2}\nonumber\\
		\le &\tfrac{n}{8}\log n+\tfrac{9}{64}n\log\log n-\left(\tfrac{3-\log 2-36\zeta'(-1)+3\log \pi-1.1108}{12}\right)n\nonumber\\
		&+\tfrac{1}{2}\log n+\left(\tfrac{1}{4}\log\pi+\tfrac{3}{4}\log 2+\tfrac{1}{2}\right)-\tfrac{1}{4n}\log 2\nonumber\\
		< &\tfrac{n}{8}\log n+\tfrac{9}{64}n\log\log n-0.8821n+\tfrac{1}{2}\log n+1.3061.\nonumber
		\end{align}
		Since $n\ge 3$, we have $\log n\le \frac{\log 3}{3}n\le 0.3663n$. Thus, we can further estimate the upper bound using $n\ge 3$ by
		$$h_{\mathrm{Fal}}(X_n)<\tfrac{n}{8}\log n+\tfrac{9}{64}n\log\log n-0.263n.$$
		The first statement of the corollary follows from these bounds.
		The second statement on the asymptotic behavior
		$$h_{\mathrm{Fal}}(X_n)=\tfrac{1}{8}n \log n+O(n\log \log n)$$
		is a direct consequence of the first statement.
	\end{proof}
	
	\section{Complex multiplication by cyclotomic fields}\label{sec_cm}
	We bound the Faltings height of an abelian variety having complex multiplication by $(\mathbb{Q}(\zeta_n),\Sigma_n)$, the canonical CM type of the $n$-th cyclotomic field $\mathbb{Q}(\zeta_n)$ for an odd integer $n\ge 3$. We refer to the introduction for the definition of the CM type $(\mathbb{Q}(\zeta_n),\Sigma_n)$ and for a brief discussion on the notion of \emph{complex multiplication}. Let us first show that the Jacobian $\mathrm{Jac}(X_n)$ always contains an abelian subvariety having complex multiplication by $(\mathbb{Q}(\zeta_n),\Sigma_n)$.
	\begin{Lem}\label{lem_subvariety}
		Let $n\ge 3$ be an odd integer and let $X_n$ be the hyperelliptic curve defined by the equation $y^2=x^n-1$. The Jacobian $J_n=\mathrm{Jac}(X_n)$ always contains an abelian subvariety $A_n\subseteq J_n$ which has complex multiplication by $(\mathbb{Q}(\zeta_n),\Sigma_n)$, the canonical CM type of $\mathbb{Q}(\zeta_n)$.
	\end{Lem}
	\begin{proof}
		With respect to the coordinates $(x,y)$ on $X_n$ induced by the equation $y^2=x^n-1$, we can define the automorphism
		$$\sigma\colon X_n\to X_n,\qquad (x,y)\mapsto (\zeta_n\cdot x,y)$$
		on $X_n$. It induces an endomorphism $\widetilde{\sigma}\colon J_n\to J_n$ on the Jacobian variety of the curve $X_n$. Since $\mathrm{End}(J_n)$ is a ring, we can consider the endomorphism $\Phi_n(\widetilde{\sigma})$ induced by the endomorphism $\widetilde{\sigma}$ and the $n$-th cyclotomic polynomial $\Phi_n(T)\in \mathbb{Z}[T]$. We define 
		$$A_n=\ker\Phi_n(\widetilde{\sigma})^0\subseteq J_n$$
		to be the abelian subvariety of $J_n$ given by the connected component of the kernel of the endomorphism $\Phi_n(\widetilde{\sigma})$.
		
		Let us prove that $A_n$ has complex multiplication by $(\mathbb{Q}(\zeta_n),\Sigma_n)$. We first show that $A_n$ has dimension $\dim A_n=\frac{\varphi(n)}{2}=\frac{[\mathbb{Q}(\zeta_n):\mathbb{Q}]}{2}$. Let $\eta_1,\dots,\eta_g$, where $g=\frac{n-1}{2}=\dim J_n$, be the basis in $H^0(J_n,\Omega_{J_n}^1)$ induced by the isomorphism $H^0(J_n,\Omega_{J_n}^1)\cong H^0(X_n,\Omega_{X_n}^1)$ and the basis $\frac{dx}{y},\frac{xdx}{y},\dots,\frac{x^{g-1}dx}{y}$ of $H^0(X_n,\Omega_{X_n}^1)$. Since $$\sigma^*\left(\frac{x^{j-1}dx}{y}\right)=\frac{(\zeta_nx)^{j-1}d(\zeta_nx)}{y}=\zeta_n^j\cdot \frac{x^{j-1}dx}{y},$$ the endomorphism $\widetilde{\sigma}$ acts on the basis $\eta_j$ by $\widetilde{\sigma}^*\eta_j=\zeta_n^j\eta_j$. Consequently, we get $$\Phi_n(\widetilde{\sigma})^*\eta_j=\Phi_n(\zeta_n^j)\eta_j.$$
		Since $\Phi_n(\zeta_n^j)=0$ if and only if $\gcd(j,n)=1$, we get
		\begin{align*}
		\dim A_n&=\dim H^0(A_n,\Omega_{A_n}^1)=\dim \ker \Phi_n(\widetilde{\sigma})^*\\
		&=\#\{j\in\{1,\dots,g\}~|~\gcd(j,n)=1\}=\frac{\varphi(n)}{2},
		\end{align*}
		where we used $g=\frac{n-1}{2}$ for the last equality.
		
		Next, we show that there is an injection of $\mathbb{Q}(\zeta_n)$ into $\mathrm{End}^0(A_n)$. Since $\mathrm{End}(A_n)$ is a ring, we have a morphism
		$$\mathbb{Z}[T]\to \mathrm{End}(A_n),\qquad P(T)\mapsto P(\widetilde{\sigma}|_{A_n}).$$
		Since $\Phi_n(\widetilde{\sigma}|_{A_n})=\Phi_n(\widetilde{\sigma})|_{A_n}=0$, this morphism factors to
		$$\mathbb{Z}[\zeta_n]=\mathbb{Z}[T]/(\Phi_n(T))\to \mathrm{End}(A_n).$$
		Since $A_n$ is a positive dimensional abelian variety, the image of this morphism contains at least $\mathbb{Z}$. It follows that the morphism $\mathbb{Z}[\zeta_n]\to\mathrm{End}(A_n)$ is injective since otherwise its image would be finite. This implies that we also get an injective map
		$$\iota\colon \mathbb{Q}(\zeta_n)=\mathbb{Z}[\zeta_n]\otimes_{\mathbb{Z}}\mathbb{Q}\xhookrightarrow{}\mathrm{End}(A_n)\otimes_{\mathbb{Z}}\mathbb{Q}=\mathrm{End}^0(A_n)$$
		after tensoring with $\mathbb{Q}$.
		
		Finally, we have to show that $T_0A_n\cong\bigoplus_{\tau\in\Sigma_n}\mathbb{C}_{\tau}$ as $\mathbb{Q}(\zeta_n)\otimes_\mathbb{Q}\mathbb{C}$-modules. We have already seen that there is an isomorphism
		$$H^0(A_n,\Omega_{A_n}^1)\cong\bigoplus_{\substack{1\le j\le g\\ \gcd(j,n)=1}}\mathbb{C}\cdot\eta_j$$
		of $\mathbb{Q}(\zeta_n)\otimes_\mathbb{Q}\mathbb{C}$-modules, where $\zeta_n$ acts on $\eta_j$ by the scalar $\zeta_n^j$. Hence, on every $\mathbb{C}\cdot \eta_j$ we get a structure of an $\mathbb{Q}(\zeta_n)$-module by sending $\zeta_n$ to $\zeta_n^j=e^{2\pi i j/n}$. Thus, it is induced by a morphism $\tau\in \mathrm{Hom}(\mathbb{Q}(\zeta_n),\mathbb{C})$ with $\mathrm{Im}\,\tau(\zeta_n)>0$. We conclude that there is an isomorphism
		$$H^0(A_n,\Omega_{A_n}^1)\cong \bigoplus_{\tau\in\Sigma_n}\mathbb{C}_{\tau}.$$		
		of $\mathbb{Q}(\zeta_n)\otimes_\mathbb{Q}\mathbb{C}$-modules. Since the restriction of holomorphic one-forms to the origin induces an isomorphism $H^0(A_n,\Omega_{A_n}^1)\xrightarrow{\sim} (T_0A_n)^\vee$, we also obtain an isomorphism $T_0A_n\cong\bigoplus_{\tau\in\Sigma_n}\mathbb{C}_{\tau}$ of $\mathbb{Q}(\zeta_n)\otimes_\mathbb{Q}\mathbb{C}$-modules by duality. This completes the proof of the lemma.
	\end{proof}
	Let us emphasize the special situation of $n=p$ being a prime.
	\begin{Rem}\label{rem_CM}
		If $n=p\ge 3$ is a prime, then it holds $$\dim \mathrm{Jac}(X_p)=\frac{p-1}{2}=\frac{\varphi(p)}{2}=\dim A_p$$ for the abelian subvariety $A_p\subseteq \mathrm{Jac}(X_p)$ in Lemma \ref{lem_subvariety}. This implies $A_p=\mathrm{Jac}(X_p)$. Thus, $\mathrm{Jac}(X_p)$ has complex multiplication by the CM type $(\mathbb{Q}(\zeta_p),\Sigma_p)$. A theorem by Colmez \cite[Theorem 0.3 (ii)]{Col93} states that the Faltings height $h_{\mathrm{Fal}}(\mathrm{Jac}(X_p))$ only depends on the CM type $(\mathbb{Q}(\zeta_p),\Sigma_p)$. Thus, Theorem \ref{thm_faltings-height-hyperelliptic} gives an explicit expression for the stable Faltings height of any abelian variety having complex multiplication by the CM type $(\mathbb{Q}(\zeta_p),\Sigma_p)$ for any prime $p\ge 3$.
	\end{Rem}
	Next, we want to bound the stable Faltings height $h_{\mathrm{Fal}}(A_n)$ for a general odd integer $n\ge 3$ to prove Corollary \ref{cor_faltings-height-cyclotomic}. Since $A_n$ in Lemma \ref{lem_subvariety} is an abelian subvariety of $\mathrm{Jac}(X_n)$, we deduce from a result by Rémond \cite{Rem22} that 	
	\begin{align}\label{equ_remond}
		h_{\mathrm{Fal}}(A_n)&\le h_{\mathrm{Fal}}(X_n)+(\dim \mathrm{Jac}(X_n)-\dim A_n)\log(\pi\sqrt{2})\\
		&=h_{\mathrm{Fal}}(X_n)+\tfrac{n-1-\varphi(n)}{2}\log(\pi\sqrt{2}).\nonumber
	\end{align}
	This bound now allows us to give the proof of Corollary \ref{cor_faltings-height-cyclotomic}.
	\begin{proof}[Proof of Corollary \ref{cor_faltings-height-cyclotomic}]
		By a theorem of Colmez \cite[Theorem 0.3 (ii)]{Col93} the stable Faltings height of an abelian variety having complex multiplication by $(\mathbb{Q}(\zeta_n),\Sigma_n)$ only depends on the CM type $(\mathbb{Q}(\zeta_n),\Sigma_n)$. Thus, we can assume that the abelian variety $A_n$ given in the statement of the corollary is obtained as the abelian subvariety $A_n\subseteq \mathrm{Jac}(X_n)$ in Lemma \ref{lem_subvariety}. In particular, we can use the bound in (\ref{equ_remond}).
		
		Let us first prove a lower bound for Euler's totient function $\varphi(n)$. Since $n\ge 3$, we get $$\varphi(n)=\prod_{\substack{p\text{ prime}\\ p\mid n}}p^{\mathrm{ord}_p(n)-1}(p-1)\ge \prod_{\substack{p\text{ prime}\\ p\mid n}}p^{\frac{\log2}{\log 3}(\mathrm{ord}_p(n)-1)}p^{\frac{\log2}{\log 3}}=n^{\frac{\log2}{\log 3}}.$$
		From the classical bound $x-1\ge \log x$ for $x>0$ we get $\frac{e^x}{e}\ge x$ by taking the exponential function and hence, $t\ge e\log t$ for all $t>1$. Since $n\ge 3$, we can conclude
		$$\varphi(n)\ge n^{\frac{\log2}{\log 3}}\ge e\cdot\tfrac{\log2}{\log3}\cdot\log n\ge 1.715\log n.$$
		Applying this to the bound in (\ref{equ_remond}), we get
		\begin{align*}
			h_{\mathrm{Fal}}(A_n)&\le h_{\mathrm{Fal}}(X_n)+\tfrac{n-1}{2}\log(\pi\sqrt{2})-\tfrac{1.715\log(\pi\sqrt{2})}{2}\log n\\
			&\le h_{\mathrm{Fal}}(X_n)+0.7457n-1.2788\log n-0.7456.
		\end{align*}
		If we combine this with the upper bound for $h_{\mathrm{Fal}}(X_n)$ in (\ref{equ_upper-bound}), we get
		\begin{align*}
			h_{\mathrm{Fal}}(A_n)&<  \tfrac{n}{8}\log n+\tfrac{9}{64}n\log\log n -0.1364n-0.7788\log n+0.5605.
		\end{align*}
		Using $n\ge 3$, we estimate
		$$h_{\mathrm{Fal}}(A_n)< \tfrac{n}{8}\log n+\tfrac{9}{64}n\log\log n-0.136n$$
		as claimed in the corollary.
	\end{proof}
	
\end{document}